\documentclass[12pt,a4paper]{article}
\usepackage[utf8]{inputenc}
\usepackage{amsmath}
\usepackage{amsfonts}
\usepackage{amssymb}
\usepackage{graphicx}
\usepackage{mathabx}

    \usepackage{pdfsync}

\author{Yves Le Jan}
\title{Markov loops topology } 
\begin{document}
\maketitle



\section{Introduction}

In the seminal work of Symanzik \cite{Symanz}, Poisson ensembles of Brownian loops were implicitly used.

 Since the work of Lawler and Werner \cite{LW} on "loop soups", these ensembles have also been the object of many investigations. Their properties can be studied in the context of rather general Markov processes, in particular Markov chains on graphs (Cf \cite{rflj}, \cite{stfl}, \cite{lejanito}, \cite{lejanbakry} )  .

The purpose of the present work is to explore their topological properties.

\section{Geodesics and loops on graphs}

We consider a finite connected graph $\mathcal{G}=(X,E))$.\\
The set of oriented edges is denoted $\overrightarrow{E}$. We also set, for any oriented edge  $\overrightarrow{e}=(e^{-},e^{+})$,  $-\overrightarrow{e}=(e^{+},e^{-})$.\\
Recall that on graphs, geodesics are defined as non backtracking paths: \\ $(x_{0},x_{1},...,x_{n})$ with $\{x_{i},x_{i+1}\}$ in $E$ and $x_{i-1}\neq x_{i+1}$.\\
Fundamental groups  $\Gamma_{x}$ are defined by geodesics from $x $ to $ x$ equipped with concatenation with erasure of backtracking subarcs. They are all isomorphic to the free group with $\vert E \vert-\vert X\vert+1$ generators, in a non canonical way. The isomorphisms, as well as a set of generators for the free group, can be defined by the choice of a spanning tree of the graph.\\
However, geodesic loops are in canonical bijection with the conjugacy classes 
of all  $\Gamma_{x}$.\\
Each loop $l$ is homotopic to a unique geodesic loop $l^g$. \\
\begin{figure}
 \includegraphics[width=6cm] {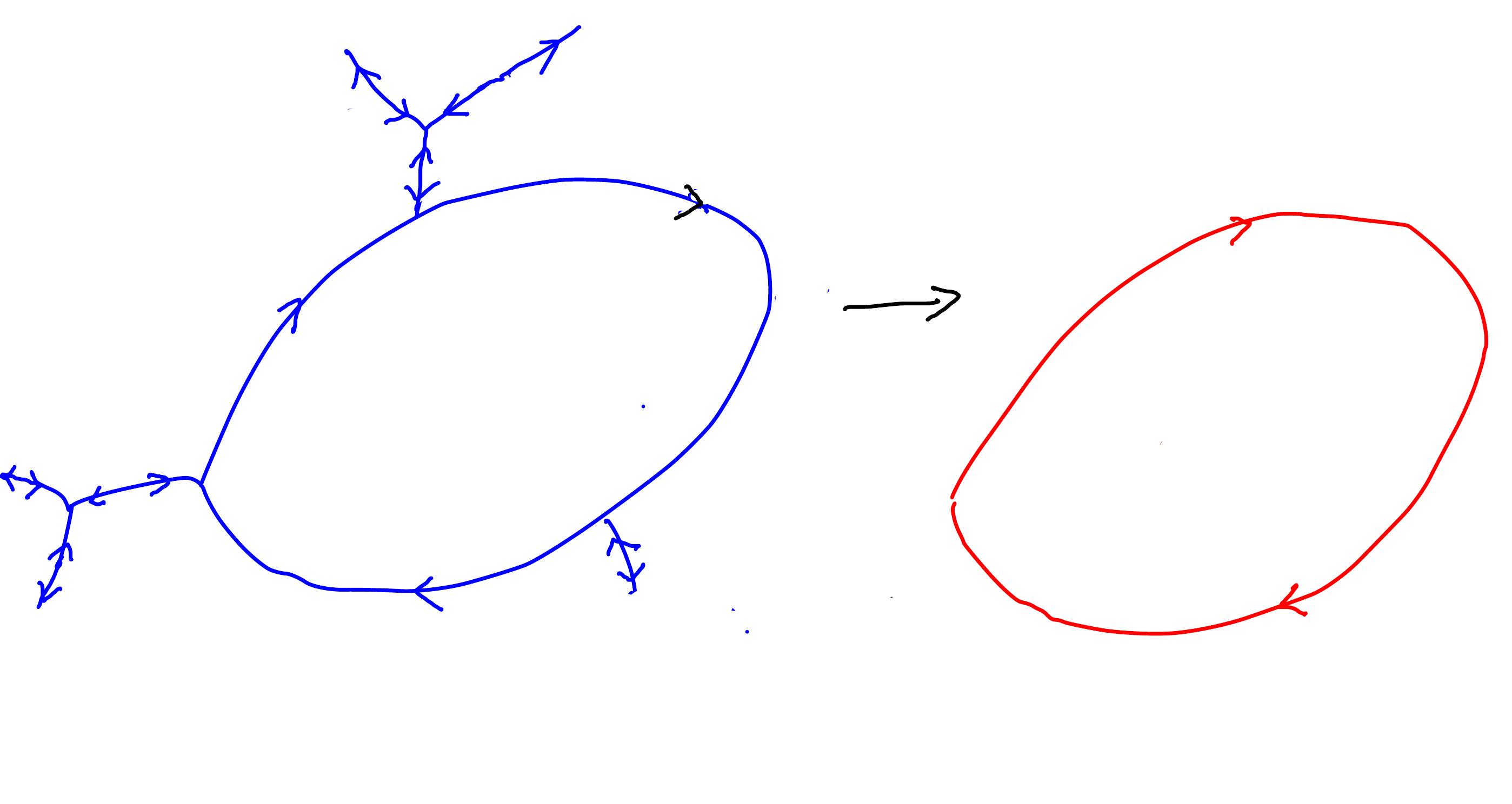}
 \caption{Loop $\rightarrow$ Geodesic Loop}
 \end{figure}
 \section{Markov loops}
 We attach a positive conductance $C_e$ to each edge $e\in E$ and a killing rate $\kappa_x$ to each vertex $x\in X$, then
define the duality measure  $\lambda_{x}=\kappa_{x}+\sum_{y}C_{x,y}$\\
and the 
 $\lambda$-symmetric transition matrix $P_y^x=\dfrac{C_{x,y}}{\lambda_x}$, $P_{\Delta}^x=\kappa_{x}{\lambda_x}$.
 The energy functional is:

$$\epsilon(f,f)=\frac{1}{2}\sum_{x,y}C_{x,y}(f(x)-f(y))^{2}+\sum_{x}\kappa_{x}f(x)^{2}$$

%
%

\medskip
We define a measure $\mu$ on (discrete time, unbased) loops:  $$\mu(l)=\frac{1}{\text{mult}(l)}\prod_{ \text{edges of } l} \left( P^{e^{-}}_{e^{+}}\right) ^{\text{mult}(l)}.$$
 Here $\text{mult}(l)$ denotes the multiplicity of the loop $l$. Note that
$$\vert \mu \vert =\mu (1)=-\log (\det(I-P))$$\\
Recall (\cite{stfl}) that this measure is induced by the restriction to non-trivial discrete loops of the measure $\sum_{x\in X}\int_{0}^{\infty}\frac{1}{t}\,\mathbb{P}_{t}^{x,x}\lambda_x dt$ defined on continuous time based loops, $\mathbb{P}_{t}^{x,x} $ being the non-normalized bridge measure defined by the transition semigroup $\exp(t[I-P])$ associated with the energy functional.\\
\medskip
A probability measure $\nu$ is defined on spanning trees (Cayley): $$\nu(T)= \frac{\prod_{ \text{edges of }T }P^{e^{-}}_{e^{+}}}{\det(I-P) }$$

%

Recall that Wilson's algorithm, based on loop erasure can be extended to provide samples $(T, \mathcal{L})$: \\ 
$T$ is a sample of $\nu$ and $\mathcal{L}$ is a sample of a Poisson point process of intensity $\mu$, i.e.  $N$ samples of $\frac{\mu}{\mu(1)}$, where $N$ is an independent Poisson variable of mean $\mu(1)$. \\
The algorithm  follows the following steps:\\
- Order $X$\\
- Run a $P$-Markov chain from  the first vertex $x_0$ to the added cemetery point $\Delta$ .\\
- Erase all loops (starting from $x_0$)  to obtain a self avoiding path $\gamma$.\\
- Restart the chain from the first point $\notin \gamma$, until the path hits $\gamma$.\\
- Erase loops and iterate until $X$ is covered. We have obtained a spanning tree $T$ and a set of based loops $\{ l_x, \: x\in X\}$.\\
- Divide each $l_x$ at its base point $x$: If $l_x$ visits $x$ $n_x$ times,\\  partition it into $\{ n_1, n_2...,n_k\}$  with probability $\frac{\prod_1^k (n_i -1)!}{n_x!}$ and define $\mathcal{L}$ to be the associated set of (unbased) loops .\\

\section{DISTRIBUTION OF GEODESIC LOOPS}

If $(x,y)$ is an edge,  let us denote $r^{x,y}$ the probability that the Markov chain starting at $y$ returns to $y$ without visiting $x$ and following a tree-contour subloop (cf Figure 1). Note that:
  $$r^{x,y}=\sum_{z\neq x}P^{y}_z P^{z}_y \sum_{n=0}^{\infty}[r^{y,z}]^n$$
Clearly, if $\gamma$ varies in the set of geodesic loops (conjugacy classes), $\vert\{l\in \mathcal{L},\, l^g=\gamma\} \vert$ are independent Poisson r.v. with mean values
$$\mu(\gamma)=\frac{1}{\text{mult}(\gamma)}(\prod_{\overrightarrow{e}\in \gamma}P_{e^+}^{e-} \rho^{e^-,e^+})^{\text{mult}(\gamma)} $$
with $\rho^{x,y}=\sum_{n=0}^{\infty}[r^{x,y}]^n$\\

\bigskip
Note that $\rho$ satisfies the relation:

$$\rho^{x,y}=1+\sum_{z\neq x}P^{y}_z P^{z}_x \rho^{x,y} \rho^{y,z}$$

\medskip
 Let us now denote $r^{x,y,k}$ the probability that the Markov chain starting at $y$ returns to $y$ for the first time in $2k$ steps following a tree-contour subloop and without visiting $x$. Set $r^{x,y}(s)=\sum r^{x,y,k}s^k$. Set $\rho^{x,y}(s)=\sum_{n=0}^{\infty}[r^{x,y}(s)]^n$\\  
 Note that:
  $$r^{x,y}(s)=s\sum_{z\neq x}P^{y}_z P^{z}_y \rho^{y,z}$$
   
\medskip
Note that $\rho^{x,y}(s)$ satisfies the relation:

$$\rho^{x,y}(s)=1+s\sum_{z\neq x}P^{y}_z P^{z}_x \rho^{x,y}(s) \rho^{y,z}(s)$$

Let us now denote $r^{x,k}$ the probability that the Markov chain starting at $x$ returns to $x$ for the first time in $2k$ steps following a tree-contour subloop. Set $r^{x}(s)=\sum r^{x,k}s^k$ Note that:
  $$r^{x}(s)=s\sum_{y}P^{x}_y P^{y}_x \rho^{x,y}(s)$$
 
  Let  denote $\rho^{x,k}$ the probability that the Markov chain starting at $x$ returns to $x$ in $2k$ steps following a tree-contour subloop. Set $\rho^{x}(s)=\sum_0 ^\infty \rho^{x,k}s^k$. Note that:  $\rho^{x}(s) = \frac{1}{1-r^{x}(s)}$ 
  and the number of loops of $\mathcal{L}$ based at $x$, homotopic to a point is a Poisson r.v. with expectation $\vert X \vert \int_0^1\frac{\rho^{x}(s)-1}{s}ds$
\medskip

If $\mathcal{G}$ is a $d$-regular graph, with $C_e=1$, $\kappa$ constant, we see that $$\rho^{x,y}(s)=\frac{(d+\kappa)^2}{2s(d-1)}(1-\sqrt{1-\frac{4s(d-1)}{(d+\kappa)^2}})$$
 $$\rho^{x}(s)=\frac{2(d-1)}{d-2+d\sqrt{1-\frac{4s(d-1)}{(d+\kappa)^2}})}$$
and recover the result of \cite{Mn} :

\medskip
$$\mu(\gamma)=\frac{1}{\text{mult}(\gamma)} \left( \frac{d+\kappa}{2(d-1)}( 1-\sqrt{1-\frac{4(d-1)}{(d+\kappa)^2}}) \right) ^{\vert \gamma \vert} $$\\
\medskip

From the expression of $\rho^x$, we deduce that number of loops homotopic to a point is a Poisson r.v. of expectation $$4\frac{d-1}{d}\vert X \vert (d( \ln(2)-\ln(b+1))+ (d-2) (\ln(b+\frac{d-2}{d})-\ln(1+\frac{d-2}{d})))$$ with $b= \sqrt{1-4\frac{d-1}{(d+\kappa)^2}}$.

\section{CONNEXIONS AND HOLONOMIES}

Recall that free groups are conjugacy separable:
Two conjugacy classes are separated by a morphism in some \textbf{finite} group G.\\  

For the  fundamental groups $\Gamma_x$ morphisms are obtained from
maps $A$, assigning to each oriented edge $ \overrightarrow{e} $ an element  $A( \overrightarrow{e} )\in G$ with $A(-\overrightarrow{e})=A(\overrightarrow{e})^{-1}.$\\
A based loop is mapped  to the product of the image by $A$ of its oriented edges and the associated loop  $l$ to the conjugacy class of this image, denoted $H_A(l)$. Moreover $H_A(l)=H_A(l^o)$
%
%
 A gauge equivalence relation between assignment maps is defined as follows:  $A_1\sim A_2$ iff  there exists $ Q$: $X\mapsto G$ such that:
 $$A_1(\overrightarrow{e})=Q(e^{+})A_1(\overrightarrow{e})Q^{-1}(e^{-})$$
Equivalence classes are $G$-connexions. They define $G$- Galois coverings of $\mathcal{G}$ (cf \cite{lejanbakry}). Obviously, holonomies depend only on the connection defined by $A$.\\
 Given a spanning tree $T$, there exists a unique  $A^T\sim A $ such that $A^T(e) =I$ for every edge $e$ of $T$.\\ For any unitary representation $\pi$ of $G$, denote $ \chi_{\pi}(C)$ the normalized trace of the image of any element in the conjugacy class $C$.\\
As $G$, $\pi$ and $A$ vary, functions $\gamma \mapsto  \chi_{\pi}(H_A(\gamma))$ span an algebra and separate geodesic loops.




 \bigskip
 Define an extended transition matrix $P^{A,\pi}$ with entries in  $X\times \{1, 2,...\dim(\pi)\}$ by $[P^{A,\pi}]^{x,i}_{y,j}=P^x_y [\pi(A(x,y))]^i _j$. Then:
 $$\sum \chi_{\pi}(H_A(l) \mu(l)=-\frac{1}{\text{dim}(\pi)}\log(\det(I-P^{A,\pi}))$$ and $\vert\{l\in \mathcal{L},\,H_A( l)= C \} \vert$ are independent Poisson r.v. with expectations:
 $$\mu(\{l,  \:H_A( l)= C\})=-\sum_{\pi \in \mathcal{R} }\overline{ \chi_{\pi}(C)}\frac{\vert C \vert}{\vert G \vert}\text{dim}(\pi) \log(\det(I-P^{A,\pi}), )$$
 $ \mathcal{R}$ denoting the set of irreducible unitary representations of $G$.

\

\section{IN THE CONTINUUM}
\
Some aspects of the theory extend to manifolds after proper rescaling and renormalisation).\\
We consider a Riemannian manifold $M$ of dimension $n$ with metric tensor $g_{i,j}$, and a potential (killing rate) $k$ on it. The energy functional is:
\[
\epsilon(f,f)=\frac{1}{2}\int g^{i,j}(x)\frac{\partial f}{\partial x_{i}%
}\frac{\partial f}{\partial x_{j}}\det(g)^{-\frac{1}{2}}dx+\int k(x)f(x)^{2}%
\det(g)^{-\frac{1}{2}}dx
\]
The heat semigroup $P_t$ is associated with the infinitesimal generator
\[
\frac{1}{2}\Delta_{x}-\kappa(x)
\]
Its kernel is denoted $p_t(x,y)$.\\

The $\sigma-$finite measure $\mu$ and the Poisson process of Brownian loops are
defined in the same way as Lawler and Werner ''loop soup'' (Cf \cite{LW}).
More precisely, $\mu=\int_{x\in X}\int_{0}^{\infty}\frac{1}{t}\mathbb{P}_{t}^{x,x} dt \det(g)^{-\frac{1}{2}}dx$ where $\mathbb{P}_{t}^{x,y}$ denotes the Brownian bridge distribution multiplied by $p_t(x,y).$

\section{HOMOTOPY CLASSES}
In this section we consider only the case of a compact surface with constant negative curvature and constant  $\kappa$ which can be represented as the quotient  $\Gamma\setminus\mathbb{H}$ of the hyperbolic plane by a discrete group of isometries  $\Gamma$.  \\
$\Gamma$ is the fundamental group of the surface and loop homotopy classes are in one-to one correspondence with closed geodesics.
If  $\gamma$ is a closed geodesic, we set
 $\vert \gamma \vert$ = length($\gamma$)\\
 It follows from an integration of the correponding term of Selberg's trace formula (cf \cite{MK}) for to the heat kernel that:
$$\mu(\{l,  \: l\text{ homotopic to } \gamma \})= \int \frac{1}{t}\, \frac{e^{-t/4}}{\sqrt{4 \pi t}}\, \frac{ \vert \gamma \vert}{\text{mult}(\gamma)}\, \frac{e^{- \vert \gamma \vert^2 /4t}}{ 2\sinh(\vert \gamma \vert/2)}\, e^{-2\kappa t}\, dt$$

Hence, setting $u=\sqrt{\dfrac{1}{4}+2\kappa}$, from the expression of the Green function of $-\Delta+(\frac{1}{4}+2\kappa)$ in $\mathbb{R}^3$, we get that:
$$\mu(\{l,  \: l\text{ homotopic to } \gamma \})= \frac{1}{4\sqrt{2\pi} \, \text{mult}(\gamma)}\dfrac{\text{K}_{u}(u\, \vert \gamma \vert)}{\sinh(\vert \gamma \vert/2)}$$.

 \section{FLAT CONNEXIONS AND HOLONOMIES}
 Given a compact Lie group $G$ with Lie algebra $\mathfrak{g}$, a $\mathfrak{g}$-valued $1$-form $A$ inducing a flat connection, a formula can be given for the distribution of the loop holonomies.
 For a smooth path $\eta$ indexed by $[0,T]$ define $X_{\eta}^A(s)$ to be the solution of the differential equation: $dX_s=A(\eta (s))X_s$.
Given a unitary representation $\pi$ of $G$, define a matrix-valued heat kernel $p_t^{A,\pi}$ with entries in  $ \{1, 2,...\dim(\pi)\}$ by $[p_t^{A,\pi}(x,y)]^{i}_{j}=\int( [\pi( X_{\gamma}^A)]^i _j)\mathbb{P}_{t}^{x,y}(d\gamma)$. \\ The definition of the multiplicative integral $X_{\gamma}^A$ on a non-smooth path could be given using Stratonovich integral but we can note instead that we can take $X_{\gamma}^ A=X_{\gamma'}^A$ if $\gamma'$ is a smooth path  close enough to $\gamma$,with time length $T(\gamma)=T(\gamma')$ and the same endpoints. Indeed, a smooth loop which is not homotopic to zero has a minimum positive diameter and if the uniform distance between two smooth paths $\gamma'$ and $\gamma"$ is small enough, we can cut them into path segments of small diameter and join the extremities of these segments by geodesics to produce a chain of loops of null holonomy from which we can deduce that $X_{\gamma'}^A=X_{\gamma"}^A$. \\

Then, if $C_m$ are disjoint central compact subsets of $G$, not containing the identity
 $\vert\{l\in \mathcal{L},\,H_A( l) \in C_m \} \vert$ are independent Poisson r.v. with expectations:
 $$\mu(\{l,  \:H_A( l) \in C_m\})=-\sum_{\pi \in \mathcal{R} } \text{dim}(\pi)^2\int _{C_m }\overline{ \chi_{\pi}(g)} dg \,\zeta_{A,\pi}'(0)$$
 
 where $dg$ denotes the normalised Haar measure on $G$, $ \mathcal{R} $ the irreducible unitary representations, $\chi_{\pi}(g)$ the normalized trace of $\pi(g)$, and $\zeta_{A,\pi}$ the meromorphic extension of the zeta function defined for $s>\frac{n}{2}$ by $$\zeta_{A,\pi}(s)=\dfrac{1}{\text{dim}(\pi)\Gamma(s)}\int_0^{\infty}t^{s-1}\text{Tr} [p_t^{A,\pi}]\,dt$$. See \cite{lejanbrow} for references and proof sktech  in the Abelian case. The proof here is similar: $\int T(l)^s\chi_{\pi}(H_A( l))-1)\mu(dl)$ is well defined and holomorphic as loops with non trivial holonomy have a minimal positive diameter, analytic continuation shows that the holomorphic functions $\zeta_{A,\pi}(s)-\zeta(s)$ and $\frac{1}{\Gamma(s)}\int T(l)^s\chi_{\pi}(H_A( l))-1)\mu(dl)$ are equal. Then \\ $\int (\chi_{\pi}(H_A( l))-1)\mu(dl)=\dfrac{d}{ds}_{ | s=0}\frac{1}{\Gamma(s)}\int T^s( \chi_{\pi}(H_A( l))-1)\mu(dl)= \zeta_{A,\pi}'(0)-\zeta'(0)$, as the reciprocal gamma function vanishes and has unit derivative in zero. Finally we conclude by Peter-Weyl theorem, noting that $\sum_{\pi \in \mathcal{R} } \text{dim}(\pi)^2\int _{C }\overline{ \chi_{\pi}(g)} dg$ vanishes as $C_m$ does not contain the identity.\\ 

\bigskip

\noindent

  D\'epartement de Math\'ematique. Universit\'e Paris-Sud.  Orsay, France.

\bigskip
   yves.lejan@math.u-psud.fr

\end{document}